\documentclass{amsart}
\usepackage[T1]{fontenc}   
\usepackage{amsfonts}
\usepackage{amsmath}
\usepackage{stmaryrd}
\usepackage{amssymb}
\usepackage{graphicx}
\usepackage{appendix}
\usepackage{subfigure}
\usepackage{url}
\usepackage{hyperref}
\usepackage{wrapfig}

\newtheorem{theorem}{Theorem}[section]

\theoremstyle{definition}
\newtheorem{definition}[theorem]{Definition}

\numberwithin{equation}{section}

\makeatletter
\newcommand*{\rom}[1]{\expandafter\@slowromancap\romannumeral #1@}
\makeatother

\begin{document}


\renewcommand{\bf}{\bfseries}
\renewcommand{\sc}{\scshape}
\vspace{0.5in}

\title[Knot Visualization Experiments for Verifiable Molecular Movies]%
{Knot Visualization Experiments for Verifiable Molecular Movies}

\author{J. Li}
\address{Department of Mathematics; University of Connecticut; Storrs, CT 06269, U.S.A.}
\email{hamlet.j@gmail.com}

\author{T. J. Peters}
\address{Department of Computer Science and Engineering, University of Connecticut, Storrs, CT 06269, U.S.A.}
\email{tpeters@engr.uconn.edu}
\thanks{The first and second authors were partially supported by NSF grants CMMI 1053077 and CNS 0923158, as well as by IBM Joint Study Agreement W1056109.  All statements here are the responsibility of the authors, not of the National Science Foundation nor of IBM.}

\author{K. E. Jordan}
\address{IBM T.J. Watson Research Center, Cambridge, MA, U.S.A.}
\email{kjordan@us.ibm.com}

\subjclass[2010]{57Q37, 57Q55, 57M25, 68R10, 92E10}

\keywords{Knot type, isotopy, low dimensional topology, topological invariance, computer animation, molecular simulation}

\begin{abstract} 
Classical topological concepts are applied to understand high performance computing simulations of molecules writhing in three dimensional space.  These simulations produce peta-bytes of floating point data, to describe 3 dimensional changes in molecular structure. A zero-th order analysis is achieved by viewing a computer animation synchronized with these changes.  The performance demands for animation of this voluminous data can become problematic, but techniques from low-dimensional topology are helpful.  The 3D molecule is reduced to a lower dimensional model of a 1-manifold, which undergoes a piecewise linear approximation for animation. An example is presented here to show how a 1-manifold and its PL approximation could come to have different embeddings as molecular writhing proceeds.  This should serve as a cautionary warning to animators to respect established sufficient conditions for topological preservation so that the movies generated will faithfully reflect the topology of the underlying model.  To obtain this result, techniques from low-dimensional topology were joined with experimental mathematics and numerical analyses.
\end{abstract}

\maketitle

\section{\bf Introduction}
\begin{figure}[h!]
\centering
      \includegraphics[height=4cm]{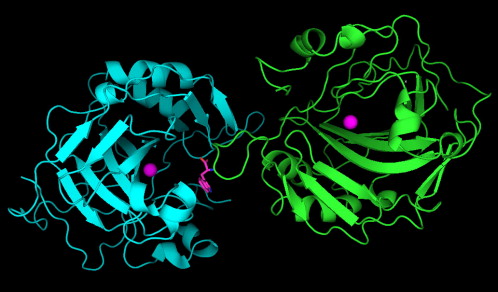}
    \caption{The backbone of a molecule \cite{temperini2008carbonic}}
    \label{fig:bb}
\end{figure}

For computations on protein molecules, a first conceptual reduction is provided by chemists through their notion of the backbone of the molecule (Figure~\ref{fig:bb}).  This is a skeletal structure which can be approximated by a 1-manifold, replacing a 3-dimensional model, to considerable computational advantage.  This 1-manifold has a computational representation as a polynomial.  For animation, this 1-manifold undergoes a piecewise linear (PL) approximation.  There are techniques to ensure that a static 1-manifold and its PL approximation are ambient isotopic.

The primary example presented in Section~\ref{sec:ecqb} synthesizes results from two primary references \cite{TJP08, JL2012} to advise caution for appropriate graphics approximations during dynamic visualization. It was motivated by two experiments conducted with our knot visualizing tool \cite{JL2012}. The first visual experiment is described in Section~\ref{sec:dr}, to observe the convergence of a B\'ezier curve to its control polygon under collinear insertion, a new spline operation defined here. The second visual experiment is described in Section~\ref{sec:ecqb}, where points were chosen and moved until the topological disparity appeared.

The term `molecular movies' has been coined to include ``\ldots cell and molecular animations \ldots '' \cite{MMsite, MMovies}.  These animations can be created by individuals or can be programmatically driven.  The shape preservation issues described here are applicable to both approaches.  The creative animation process for individuals is supported by software \cite{MMaya}, which include shape changing techniques called `morphing'.

This preservation of topology is one property being explored under {\em verifiable visualizations} \cite{KiSi08}: ``That will consider both the errors of the individual visualization component within the scientific pipeline and the interaction between and interpretation of the accumulated errors generated in the computational pipeline, including the visualization component.''

\subsection{Related Work}
There is contemporary interest \cite{Amenta2003, Chazal2005, JL-isoconvthm, JL-bez-iso} to preserve topological characteristics such as homeomorphism and ambient isotopy between an initial geometric model and its approximation. The models here are 1-manifolds, represented by class of polynomials known as B\'ezier curves (Please see Section~\ref{sec:dr}), where these B\'ezier curves are defined by finitely many \emph{control points} which become input to a standard PL approximation algorithm known as subdivision. This subdivision algorithm is widely used in graphics and animation, so that a smooth curve is rendered from its PL control polygon. Sufficient conditions for a homeomorphism between a B\'ezier curve and its control polygon have been studied \cite{TJP08, JL-ang-conv, M.Neagu_E.Calcoen_B.Lacolle2000}, while topological differences have also been shown \cite{Bisceglio, JL2012, Piegl}.   

\begin{figure}[h!]   
\centering
       \subfigure[ Simple $\mathcal{L}$ \& Self-intersecting $\mathcal{C}$]
    {   \includegraphics[height=3cm]{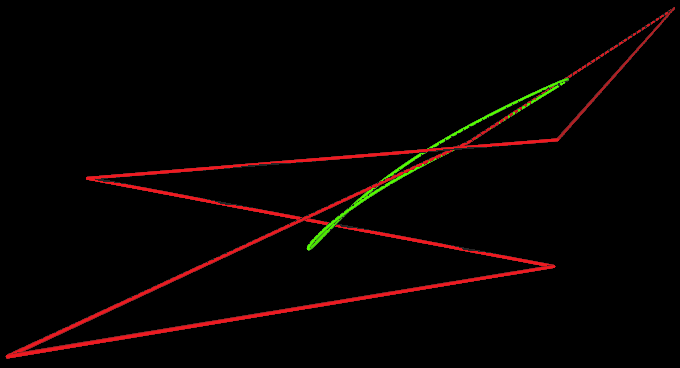}   \label{fig:eq}   }
               \subfigure[Zoomed-in intersection]
   {   \includegraphics[height=3cm]{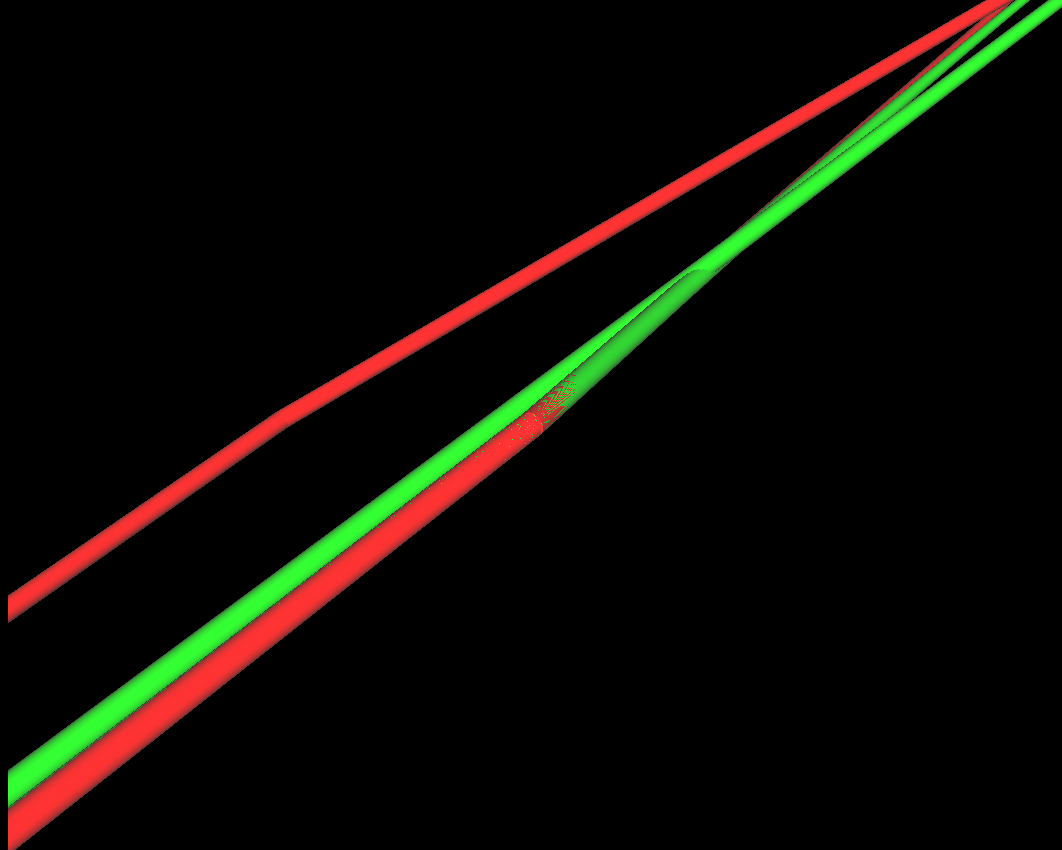}  \label{fig:eqz} }

    \caption{Equilateral, simple $\mathcal{L}$ with self-intersecting $\mathcal{C}$}    \label{fig:knots1}    
 \end{figure}

The paper \cite{JL2012} presents techniques to create a class of examples -- where a B\'ezier curve and its control polygon are not homeomorphic, as the B\'ezier curve is self-intersecting while the control polygon is simple, closed and equilateral, as shown in Figure~\ref{fig:knots1}.  These were all static examples. 

The work \cite{TJP08} also established sufficient conditions for preservation of knot type during dynamic visualization of ongoing molecular simulations, as the knot type of the molecule is an important consideration for these biochemists.  However, it relied upon a new approximation for each iteration of the perturbed curve.  These approximation algorithms are more intensive than just continuing to perturb the initial PL approximation, as is examined here.  

Upper bounds were given \cite{TJP08} on the perturbations of the vertices or a broad class of PL approximations, so as to retain the same knot type of an underlying curve model representing the writhing molecule.  It is easy to understand that drastic movements of the vertices of the PL model could create a different knot type in the graphics versus the perturbing curve, but it is shown that even small vertex changes beyond these established upper bounds could lead to incorrect knotting in the graphics display.  Appreciating this subtlety is important to assess the pragmatic trade-offs of maintaining one PL approximation and perturbing it over multiple frames versus the more performance intensive approach of re-approximating the B\'ezier curve for each frame.  

\subsection{Background}
\begin{figure}[h!]   
\centering
       \subfigure[Before Perturbations]
    {   \includegraphics[height=4.5cm]{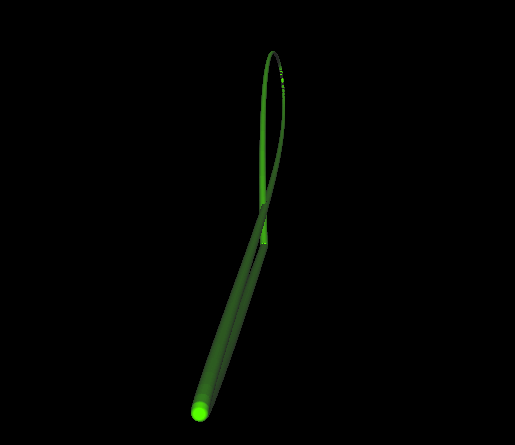}   \label{fig:before}   }
               \subfigure[Two Perturbations]
   {   \includegraphics[height=4.5cm]{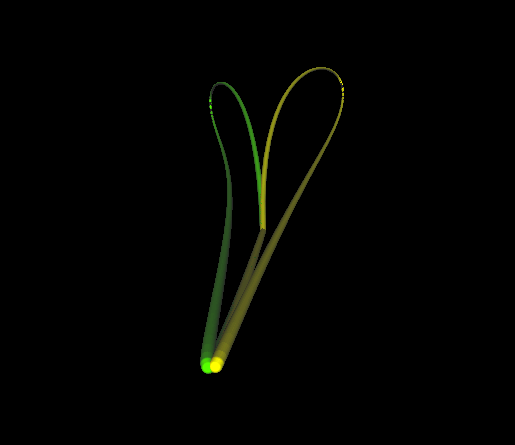}  \label{fig:after} }
    \caption{Perturbations}    \label{fig:ptc}    
 \end{figure}

The strategy chosen is to perturb some control points, while creating a self-intersection of the B\'ezier curve, as described more fully in Section~\ref{sec:id}. An immediate difficulty encountered is illustrated in Figure~\ref{fig:ptc}. Figure~\ref{fig:before} shows the B\'ezier curve from Figure~\ref{fig:knots1}. Figure~\ref{fig:after} shows two resultant B\'ezier curves after a control point is perturbed to two different positions. Since both left and right images of Figure~\ref{fig:after} appear to be unknotted, it is possible that these unknots can be transformed into each other without ever passing through an intersection. This motivates us to create an example of perturbation changing the knot type of the underlying curve, which will be presented below.  

\subsection{Ideas}\label{sec:id}
\begin{figure}[h!]
\centering
      \includegraphics[height=2.5cm]{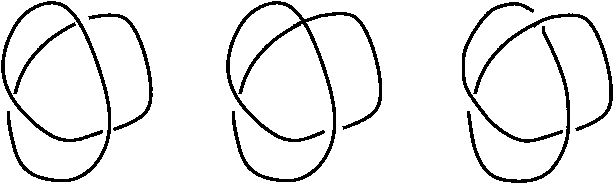}
    \caption{Necessary self-intersection}
    \label{fig:ut}
\end{figure}

The key here is to prove the existence of self-intersection by demonstrating a change in knot type. Consider Figure~\ref{fig:ut}, having the unknot indicated on the left, as the initial configuration. If the unknot is perturbed continuously to create the trefoil knot shown on the right of Figure~\ref{fig:ut}, then a self-intersection must necessarily have occurred, as shown in the middle image of Figure~\ref{fig:ut}. The sufficient condition for the self-intersection is a change of isotopy class between the original and final knot.

We generate a control polygon with a B\'ezier curve of high degree, by a new technique: collinear insertion\footnote{This terminology is used in other fields \cite{louis1997cloning}, but with a very different meaning that should not cause confusion within this mathematical context.} (described in detail in Section~\ref{sec:dr}). Collinear insertion was developed by experimental mathematical methods, as assisted by a visualization tool, which is described elsewhere \cite{JL2012}. The experiments showed that perturbing one vertex of the subdivided control polygon typically left all other vertices fixed, so that embedding changes were not likely to occur.   However, it was observed that moving the end point of a segment containing many collinear control points would cause more global changes.  Collinear insertion was conceived to provide a close approximation between a B\'ezier curve and its control polygon (as happens after sufficient subdivision) while having the \emph{additional} property of many collinear control points.

With the carefully constructed example, we can perturb a vertex of the control polygon such that the control polygon remains simple during the perturbation, while the knotting structure of the B\'ezier curve changes (i.e. an unknot is perturbed to a nontrivial knot). Hence, the corresponding B\'ezier curve must become self-intersecting during this perturbation. The change of the knotting structure relies upon a numerical analysis, where the numerical tolerances can be easily adjusted to be within a single pixel resolution, merely by scaling the dimensions of the input according to the screen size, as clearly sufficient for computer animation. 

Admittedly, this is only a first order analysis but is likely sufficient for animation, but, more importantly, indicates the opportunity for more sophisticated numerical analyses, if deemed important for other applications.  For instance, the specific code implementations of the simplex method \cite{Lagarias1998} and Horner's method \cite{cormen2001introduction} invoked in Section~\ref{sec:pu} could be subjected to a rigorous forward error analysis \cite{Scheid} to assess robustness, but such investigations are beyond the scope of this paper.  Related error bounds on spline curves have appeared as interval splines \cite{Tuohy1997791}, which could afford alternative opportunities for other, more refined analyses. While the numerical details could vary, the topological perspective presented would adapt to those subtleties. 

\subsection{Techniques: Collinear Insertion and Perturbations} Collinear insertion is similar to degree elevation \cite{G.Farin1990} in that both methods produce higher degree B\'ezier curves. However, collinear insertion preserves the perimeter of the control polygon while changing the B\'ezier curve (see Figures~\ref{fig:4c} \&~\ref{fig:rasd}), whereas degree elevation changes the control polygon but preserves the B\'ezier curve. 

Motivated by applications for dynamic visualization in a high performance computing environment, the preservation of topological integrity during perturbations has been investigated \cite{L.-E.Andersson2000, TJP08}. For example, when visualizing a molecule twisting and writhing under local chemical and kinetic changes, it is crucial that topological artifacts are not introduced by the visual approximations \cite{TJP08}. In particular, a perturbation on a control point of a B\'ezier curve changes both the control polygon and the B\'ezier curve. The paper \cite{L.-E.Andersson2000} gives an upper bound on the perturbation of vertices of polyhedra to retain ambient isotopic equivalence.   The paper \cite[Proposition 5.2]{TJP08} provides a sufficient condition to ensure that both a control polygon and its B\'ezier curve remain ambient isotopic during the perturbation. 

We present an example, consistent with previous upper bounds \cite{L.-E.Andersson2000} on preservation of ambient isotopic equivalence, where a slight perturbation of a control point results in an ambient isotopic control polygon while changing the knot type of the B\'ezier curve. 

\subsection{Outline}
In Section~\ref{sec:dr}, we describe the collinear insertion technique. In Section~\ref{sec:ecqb}, we provide images of the simple, closed control polygons and B\'ezier curves constructed by the collinear insertion and perturbation. In Sections~\ref{sec:puss} and~\ref{ssec:pk}, we establish the knot types of the original and perturbed B\'ezier curves. In Section~\ref{ssec:tscp}, we verify the non-self-intersection of the control polygon.

\section{\bf Collinear Insertion for B\'ezier Curves}\label{sec:dr}
\begin{definition} \textup{\cite{Piegl}}
\label{def:bez}
The parameterized \textbf{B\'ezier curve}, denoted as $\mathcal{B}(t)$, of degree $n$ with control points $P_m \in\mathbb{R}^3$ is defined by
$$\mathcal{B}(t) = \sum_{m=0}^{n} \binom{n}{m} t^m (1-t)^{n-m} P_m, t\in[0,1],
$$
The $PL$ curve determined by the points $\{P_0,P_1,\ldots,P_n\}$ is called the {\em control polygon}. If $P_0=P_n$, then the polygon and the curve are {\em closed}. Otherwise, they are {\em open}.
\end{definition}

During the study of the existence of an equilateral control polygon with a self-intersecting B\'ezier curve, we tried to ``slightly'' perturb a vertex so that the control polygon would remain simple, but simultaneously the B\'ezier curves would differ ``dramatically'' by a change of  knot type. However, for many examples, if a vertex is ``slightly'' perturbed, then  the B\'ezier curve is ``slightly'' changed, i.e. the knot type remains the same.  

One may use the de Casteljau algorithm to create a new control polygon which is closer \cite{Nairn-Peters-Lutterkort1999} to the B\'ezier curve in Hausdorff distance, but this also shortens the length of edge of the control polygon, such that the perturbation of a vertex only effects a small region of the control polygon. Consequently, the perturbation may not change the B\'ezier curve sufficiently to change the knot type.

\begin{figure}[h!]
\centering
      \includegraphics[height=4cm]{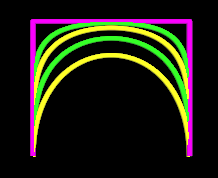}
    \caption{Successive Collinear Insertions}
    \label{fig:4c}
 \end{figure}


A fundamental intuition about B\'ezier curves is that the curve tends towards the control points.  It follows that adding more collinear control points should draw the curve closer to those edges. We were able to corroborate that visually through experiments, as shown in Figure~\ref{fig:4c}, where the sequence of smooth curves are obtained from successive collinear insertions. The innermost curve is defined by four control points, whereas the curve closest to the control polygon has 25 control points, with the intermediate ones defined by 7 and 13 control points, respectively. That experimental evidence was sufficient for creation of this example, while leaving a formal convergence proof as a subject for further investigation. 

The illustrative example of Figure~\ref{fig:4c} is of an open curve. The rest of this section treats closed curves. Suppose that there is an initial set of control points, denoted as $\{v_0,v_1,\ldots, v_{n-1}, v_0\}$, defining a closed curve of degree $n$. (Note that the vertices $\{v_0,v_1,\ldots, v_{n-1}\}$ with $v_0 \neq v_{n-1}$ would define an `open' curve of degree $n-1$. ) Now we describe the procedure to raise the degree by adding collinear control points. We designate this as \textit{collinear insertion}. 

\begin{enumerate}
\item Take the midpoints\footnote{For ease of exposition, the collinear points inserted are all midpoints. However, it is clear that other schemes for insertion of collinear points are possible.  This is similar to the typical exposition of the de Castlejau algorithm \cite{G.Farin1990}, where the midpoints are used to explain the concept, but other intermediate points are equally possible.} along each line segment. 

\item Add these midpoints to the initial set of control points to form a new set of control points: 
$$\{v_0,\frac{v_0+v_1}{2}, v_1,\ldots, v_{n-1}, \frac{v_{n-1}+v_0}{2}, v_0\}.$$ 
Consequently, the perimeter of the control polygon remains the same, but the B\'ezier curve changes. (Retaining this original perimeter, even as more control points are added, supports an easy method to perturb many control points simultaneously while only perturbing the original vertices.) The degree is raised from $n$ to $2*n$.

\item If the process is repeated $k$ times, then we obtain a set of $2^k*n$ control points, which determines a closed B\'ezier curve of degree $2^k*n$.
\end{enumerate}

\begin{figure}[h!]
\centering
    \subfigure[Initial curves]
    {
   \includegraphics[height=4.5cm]{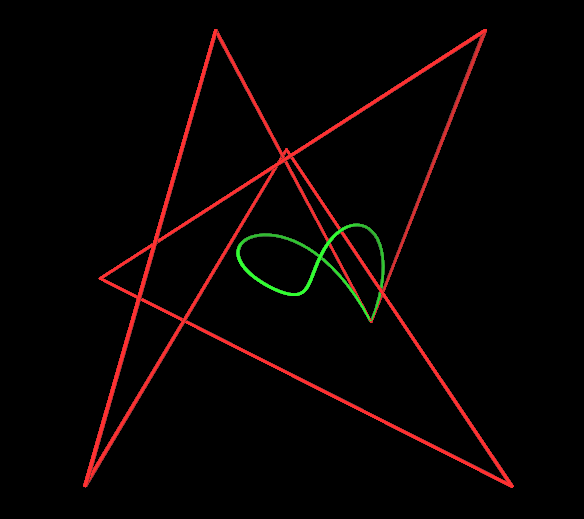}
    \label{fig:inc}
    }
    \subfigure[Higher degree curves]
    {
   \includegraphics[height=4.5cm]{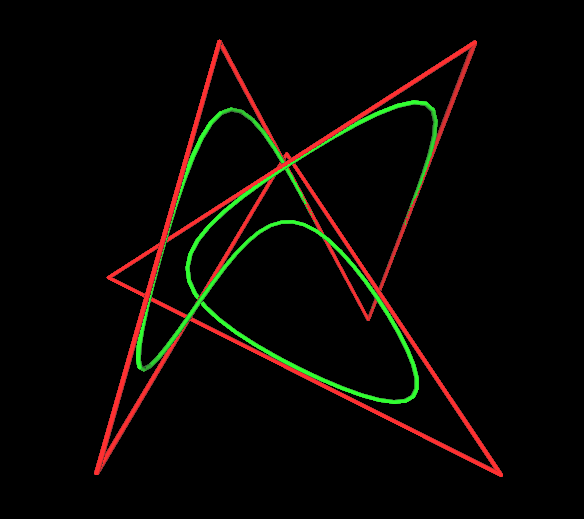}
   \label{fig:hd}
    }
    \caption{Raising the degree}
    \label{fig:rasd}
\end{figure}

Figure~\ref{fig:inc}  shows a closed control polygon with a closed B\'ezier curve of degree $7$ determined by $7$ vertices, and Figure~\ref{fig:hd} shows the control polygon with a B\'ezier curve of degree $112$ determined by $112$ control points, constructed by the collinear insertion, which is the focus of the rest of this paper.

\section{\bf The Defining Data for the Example}\label{sec:ecqb}

We start with $7$ control points\footnote{At least four decimal digits are provided here. If higher accuracy is desired, more digits may be chosen.} listed below to determine a closed control polygon. 
{\small $$(1.9817,   -1.7646,   -4.5897), (-1.3841185, 4.6825505, 0.913541),$$
$$(-3.2983075, -4.0566825, 2.686189), (-0.1232995, 2.768254, -2.463584),$$
$$(3.9079915, -4.533357, 1.2263705), (-3.935983, -0.438272, -0.983365),$$
$$(3.218174, 4.296123, 2.1124595).$$}

\begin{figure}[h!]
\centering
    \subfigure[Initial polygon and curve]
    {
   \includegraphics[height=4.5cm]{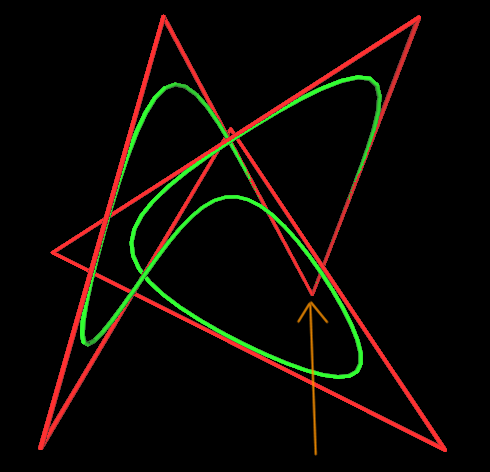}
   \label{fig:s12-a}
    }
    \subfigure[Perturbed polygon and curve]
    {
   \includegraphics[height=4.5cm]{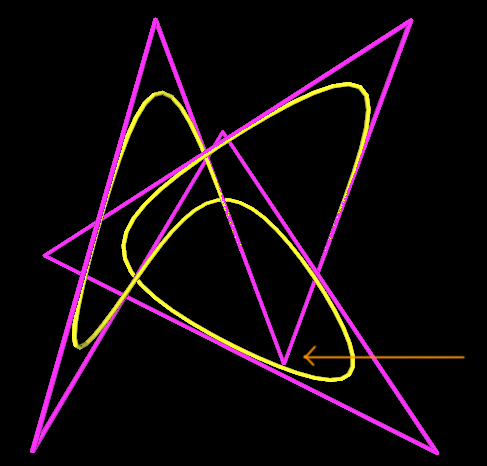}
    \label{fig:s12-b}
    }
    \caption{Perturb $v_0$ to $v'_0$}
    \label{fig:s12}
\end{figure}

We then use the collinear insertion to insert midpoints. Repeating this process four times, we obtain $112$ control points. The corresponding control polygon and B\'ezier curve of degree $112$ are shown in Figure~\ref{fig:s12-a} which repeats Figure~\ref{fig:hd} for a visual comparison to Figure~\ref{fig:s12-b}, which is obtained by perturbing one of the original vertices, as described below. Denote the control polygon and the B\'ezier curve in Figure~\ref{fig:s12-a} by $\mathcal{P}_1$ and $\mathcal{B}_1$. 

Denote the vertices of $\mathcal{P}_1$ (the original control points) as $v_0, v_1, \ldots, v_6$. Note that the vertices adjacent to $v_0$ are $v_1$ and $v_6$. We fix $v_1$ and $v_6$ while moving the line segments $\overline{v_0v_1}$ and $\overline{v_0v_6}$, by perturbing $v_0$ from $(1.9817,   -1.7646,   -4.5897)$ to $(1.3076,   -3.3320,   -2.5072)$, denoted by $v'_0$, as shown in Figure~\ref{fig:s12-b}. Denote the resultant control polygon and B\'ezier curve as $\mathcal{P}_2$ and $\mathcal{B}_2$.

The perturbation of $v_0$ to $v'_0$ is an affine transform, performed by the elementary graphics operation of rotation of a point about an arbitrary line \cite{Theorharis}.  The line of rotation is formed by the points $v_1$ and $v_6$.  This creates a new control polygon, where only two segments have been altered, those containing $v'_0$, which are $\overline{v_1v'_0}$ and $\overline{v_6v'_0}$.

When $v_0$ is perturbed to $v'_0$, the control polygon remains simple, but the unknot $\mathcal{B}_1$ is changed to a nontrivial knot $\mathcal{B}_2$ (See Section~\ref{ssec:pk}). The different knotting structures between $\mathcal{B}_1$ and $\mathcal{B}_2$ imply that there must be an instant during the perturbation, when the B\'ezier curve is self-intersecting. This is an example of a perturbation which preserves the ambient isotopy of the control polygon but changes the knot type of the B\'ezier curve. 

\section{\bf Topological Differences}\label{sec:pu}

The following three Subsections show that
\begin{itemize}
\item the initial B\'ezier curve is unknotted,
\item the perturbed B\'ezier curve is a nontrivial knot, indicating that a self-intersection of the B\'ezier curve \emph{must} have occurred during this perturbation, and
\item the control polygon remains simple throughout the perturbation.
\end{itemize}

We prove the knottedness and unknottedness under the sub-pixel criterion, that is, we accept results generated by numerical methods if any numerical tolerances can be adjusted to provide sub-pixel resolution, as sufficient for this context. 


\subsection{Initial Unknottedness}\label{sec:puss}
\label{ssec:initun}

\begin{figure}[h!]
\centering
    \subfigure[Initial curve]
    {
   \includegraphics[height=4.5cm]{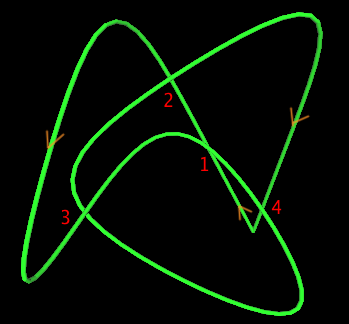}
   \label{fig:s12-ac}
    }
    \subfigure[Perturbed curve]
    {
   \includegraphics[height=4.5cm]{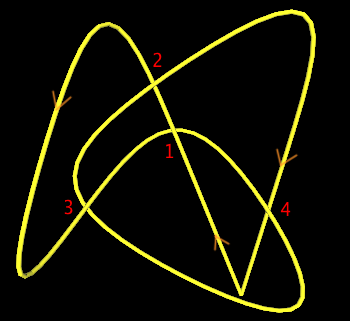}
   \label{fig:s12-bc}
    }
    \caption{UnKnot and Knot}
    \label{fig:s12c}
\end{figure}
For the unknottedness of $\mathcal{B}_1$, consider the projection shown in Figure~\ref{fig:s12-ac}, obtained by taking each z-coordinate to be zero, that is, the projection onto the xy-plane.  Denote the projected B\'ezier curve by $C2d(t)$. Define a function $fnS(t_1,t_2): [0,1] \times [0,1] \rightarrow R$ by 
$$fnS(t_1,t_2)=||C2d(t_1)-C2d(t_2)||.$$ If $fnS(t_1,t_2)=0$ and $t_1 \neq t_2$, then we have $C2d(t_1)=C2d(t_2)$ and $t_1 \neq t_2$, which implies that the projected B\'ezier curve is self-intersecting, and the self-intersection has $t_1$ and $t_2$ as parameters. In other words, we seek the roots of $fnS(t_1,t_2)$ with $t_1 \neq t_2$ to determine the self-intersections. Note that $fnS(t_1,t_2)$ is a polynomial of degree $112$. We use the simplex search method \cite{Lagarias1998} implemented by Mathlab function `fminsearch' to find the roots of this high degree polynomial. 

The pairs of parameters where the self-intersections occur are listed below, which correspond to intersections labeled as `1, 2, 3' and `4' in Figure~\ref{fig:s12-ac}:\\
{\small $$[0.0488, 0.4614], [0.0861, 0.7918], [0.3473, 0.6931], [0.5126, 0.9915].$$}

Using these parameter values to evaluate\footnote{This evaluation is performed by Matlab, where Horner's method \cite{cormen2001introduction} is invoked.} the original B\'ezier curve in 3D (Definition~\ref{def:bez}),  we get the pairs of points corresponding to the above parameter values:\\
{\small $$[(0.8309,    0.4397,   -2.7081), (0.8308,    0.4397,   -1.2072)];$$
$$[(-0.0435,    2.0929,   -1.2807), (-0.0435,    2.0928,    0.6747)];$$
$$[(-1.8672,   -1.0031,    0.4288), (-1.8672,   -1.0032,   -0.3359)];$$
$$[(2.0548,   -1.4062,   -0.3480), (2.0548,   -1.4061,   -4.1932)].$$}

Note that the orientation of $\mathcal{B}_1$ is shown by arrows in Figure~\ref{fig:s12-ac}. By comparing the z-coordinates of points listed above, we find that `1' and `2' are under crossings, while `3' and `4' are over crossings. It follows that $\mathcal{B}_1$ is simple, within the numerical limits imposed by this finite precision data and these floating point computations. 

Observe specifically that the sub-curve from `3'  to `4' is over the other portions. So using the published notation \cite{Livingston1993}, a Reidemeister move of Type 2b will eliminate the crossings `1, 3' and `4', by pulling the sub-curve from `3'  to `4' downwards. Consequently the resultant diagram after the Reidemeister move has only one crossing, which implies $\mathcal{B}_1$ is unknotted. 

\subsection{Nontrivial Knot Occurrence}\label{ssec:pk}
For the nontrivial knottedness of $\mathcal{B}_2$, we start from a projection obtained by taking each z-coordinate to be zero, as shown in Figure~\ref{fig:s12-bc}. 

Again, using the MathLab function 'fminsearch', we find the following pairs of parameters where the projected curve is self-intersecting:
$$[0.0761,    0.4240], [0.0928,    0.7830], [0.3473,    0.6931], [0.5039,    0.9575].$$

Evaluating (again, using Horner's method) the original B\'ezier curve in 3D (Definition~\ref{def:bez}) at these parameters produces the corresponding pairs of points:
{\small $$[(-0.1265,    0.9308,   -0.6850), (-0.1266,    0.9308,   -1.2854)];$$
$$[(-0.4389,    1.8160,   -0.2901), (-0.4389,    1.8159,    0.5159)];$$
$$[(-1.8672,   -1.0032,    0.4289), (-1.8671,   -1.0032,   -0.3358)];$$
$$[(1.8761,   -1.0622,   -0.5163), (1.8761,   -1.0623,   -1.1327)].$$}

\begin{figure}[h!]
\centering
      \includegraphics[height=5cm]{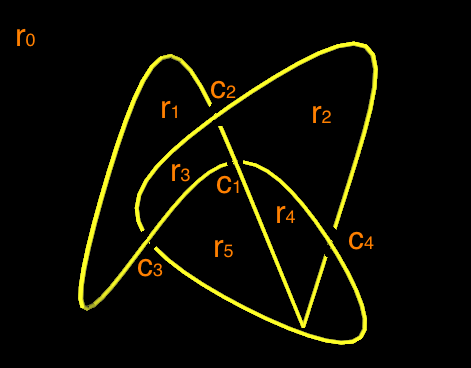}
    \caption{Diagram}
    \label{fig:diag}
\end{figure}

The under and over crossings follow from comparing the z-coordinates in each pair, as shown by the diagram in Figure~\ref{fig:diag}.

In the diagram, the crossings are labeled as $c_1, \ldots ,c_4$ and regions are $r_0, \ldots, r_5$. According to \cite{alexander1928topological}, we obtain the following corresponding matrix of the diagram:
$$
\begin{bmatrix}
  0 & 0 & 1  & -1  & -t & t\\
  -1 & t & 1  & -t & 0 & 0\\
  -t & t & 0 & -1 & 0 & 1\\
  t  & 0 & 1 & 0 & -1 & t
 \end{bmatrix}
$$

Deleting any two consecutive columns of the above matrix, computing the determinant and normalizing the polynomial, we get the Alexander polynomial $1-3t+t^2$, which implies the knot is nontrivial  \cite{alexander1928topological}. 

\subsection{Control Polygon Remains Simple}\label{ssec:tscp}
Let $\alpha$ be the angle of rotation from $v_0$ to $v'_0$. Let $v(\theta)$ be the point between $v_0$ and $v'_0$ during the rotation, where $\theta \in [0, \alpha]$ and it indicates the angle of rotation from $v_0$ to $v(\theta)$. This rotation function is the standard rotation matrix from an introductory graphics course \cite{Theorharis}. Now, to show the non-self-intersection of the control polygon during the perturbation, it suffices to show that the line segments $\overline{v(\theta)v_1}$ and $\overline{v(\theta)v_6}$ do not intersect with other line segments of the control polygon, for all $\theta \in [0,\alpha]$.  

For each $i=1,\cdots,5$, consider $\overline{v(\theta)v_1}$ and $\overline{v_iv_{i+1}}$. Parametrize the line segments by:
\[
\left\{ 
  \begin{array}{l l}
   \overline{v(\theta)v_1}: v(\theta)+(v_1-v(\theta))s, \ s\in[0,1] \\
    \overline{v_iv_{i+1}}: v_i+(v_{i+1}-v_i)t, \ t\in[0,1]\\
  \end{array} \right.
\]

And then consider the equation 
$$v(\theta)+(v_1-v(\theta))s=v_i+(v_{i+1}-v_i)t.$$
It has no solution for  $\theta,s$ and $t$ with $\theta \in [0,\alpha]$ and $s, t \in [0,1]$, so $\overline{v(\theta)v_1}$ does not intersect $\overline{v_iv_{i+1}}$. 

We can similarly verify that $\overline{v(\theta)v_6}$ does not intersect $\overline{v_iv_{i+1}}$ for each $i=1,\cdots,5$.

\section{\bf Conclusion} 

In this work, classical topology is complemented by experimental mathematics and numerical analysis to provide insight into pathologies that could occur when creating synchronous visualizations of molecular simulations. A representative example is presented where a small data perturbation produces different knot types between a polynomial curve and its PL approximation used for graphics rendering. The disparity exhibited serves as a cautionary note to the application community for rigorous attention to graphics approximations during dynamic visualization. The analyses assume a sub-pixel criterion, but, more importantly, the methods (including collinear insertion) are likely to serve as a catalyst for domain scientists to investigate the relationship between smooth and PL knotting structures. 

\section{\bf Acknowledgement} 
The authors acknowledge, with appreciation, an insightful question posed by an anonymous referee of the authors' previous paper [11], where that question led to the results presented here.

\bibliographystyle{plain}
\bibliography{ji-tjp-biblio}

\end{document}